\documentclass[12pt,twoside]{article}
\usepackage{graphicx}
\usepackage{blindtext}
\usepackage[T1]{fontenc}
\usepackage[utf8]{inputenc}
\usepackage[top=.8in, bottom=.8in, left=.7in, right=.7in]{geometry}
\usepackage{amsfonts,amsmath,mathtools,amsthm, amscd, amssymb, color,enumerate}
\usepackage{cite}
\usepackage{graphics}
\def\bpsp{\begin{pspicture}}
\def\epsp{\end{pspicture}}

\newtheorem{theorem}{Theorem}[section]
\newtheorem{remark}[theorem]{Remark}
\newtheorem{example}[theorem]{Example}
\newtheorem{lemma}[theorem]{Lemma}
\newtheorem{corollary}[theorem]{Corollary}
\newtheorem{definition}[theorem]{Definition}
\newtheorem{proposition}[theorem]{Proposition}
\newtheorem{note}{Note}
\newtheorem{case}{Case}
\newtheorem{conjecture}{Conjecture}
\newtheorem{question}{Question}

\newcommand{\bea}{\begin{eqnarray}}
\newcommand{\eea}{\end{eqnarray}}
\newcommand{\beq}{\begin{eqnarray*}}
\newcommand{\eeq}{\end{eqnarray*}}

\def\m4{\mbox{\rm ~(mod $4$)}}

\def \bd{\begin{definition}}
\def \ed{\end{definition}}
\def \bqu{\begin{question}}
\def \equ{\end{question}}
\def \bcc{\begin{conjecture}}
\def \ecc{\end{conjecture}}
\def \bt{\begin{theorem}}
\def \et{\end{theorem}}
\def \bl{\begin{lemma}}
\def \el{\end{lemma}}
\def \bc{\begin{corollary}}
\def \ec{\end{corollary}}
\def \be{\begin{equation}}
\def \ee{\end{equation}}
\def \ben{\begin{enumerate}}
\def \een{\end{enumerate}}
\def \ba{\begin{array}}
\def \ea{\end{array}}
\def \bp{\begin{proposition}}
\def \ep{\end{proposition}}
\def \bx{\begin{example}}
\def \ex{\end{example}}
\def \br{\begin{remark}}
\def \er{\end{remark}}
\def \bdsc{\begin{description}}
\def \edsc{\end{description}}

\def \bn{\begin{case}}
\def \en{\end{case}}
\def \bnt{\begin{note}}
\def \ent{\end{note}}
\def\1{1\!\!1}

\def\mm2{\mbox{\rm ~(mod $2$)}}
\def\m4{\mbox{\rm ~(mod $4$)}}

\def\qed{\nolinebreak\hfill\rule{.2cm}{.2cm}\par\addvspace{.5cm}}

\def\m{\mu}

\def\1{\textbf{1}}
\def\0{\textbf{0}}

\begin{document}
\title{Spectra  of s-neighbourhood corona of two  signed graphs  }
\author{ Tahir Shamsher$^{a}$, Mir Riyaz ul Rashid$^{b}$ and S. Pirzada$^{c}$ \\
$^{a,b,c}${\em Department of Mathematics, University of Kashmir, Srinagar, Kashmir, India}\\
 \texttt{tahir.maths.uok@gmail.com}; \texttt{mirriyaz4097@gmail.com}\\
 \texttt{pirzadasd@kashmiruniversity.ac.in}}
\date{}

\pagestyle{myheadings} \markboth{Tahir, Riyaz, Pirzada }{Spectra  of s-neighbourhood corona of two  signed graphs   }
\maketitle
\vskip 5mm
\noindent{\footnotesize \bf Abstract.} A signed graph $S=(G, \sigma)$ is a  pair  in which $G$ is an underlying graph and $\sigma$ is a function from the edge set to $\{\pm1\}$.  For signed graphs $S_{1}$ and $S_{2}$  on $n_{1}$ and $n_{2}$ vertices, respectively,  the signed neighbourhood corona $S_{1} \star_s S_{2}$ (in short s-neighbourhood corona)  of $S_{1}$ and $S_{2}$ is the signed graph obtained by taking one copy of $S_{1}$ and $n_{1}$ copies of $S_{2}$  and joining every neighbour of the $i$th vertex of $S_{1}$  with the same sign as the sign of incident edge to every vertex in the $i$th copy of $S_{2}$. In this paper, we investigate the adjacency, Laplacian and net Laplacian spectrum of $S_{1} \star_s S_{2}$ in terms of the corresponding spectrum of $ S_{1}$ and $ S_{2}$. We determine
$(i)$ the adjacency spectrum of $S_{1} \star_s S_{2}$ for arbitrary $S_{1} $ and net regular $ S_{2}$, $(ii)$ the Laplacian spectrum  for regular $S_{1} $ and regular and net regular $ S_{2}$ and $(iii)$ the net Laplacian spectrum  for net regular $S_{1} $ and arbitrary $ S_{2}$. As a consequence, we obtain the signed graphs with $4$ and $5$ distinct adjacency, Laplacian and net Laplacian eigenvalues. Finally, we show that the signed neighbourhood corona of two signed graphs is not determined by its adjacency (resp., Laplacian, net Laplacian) spectrum.

\vskip 3mm



\section{Introduction}\label{sec1}
 A signed graph is defined to be a pair $S=(G,\sigma)$, with $G=(V_G,E_G)$ as the underlying graph and $\sigma: E_G \rightarrow \{-1, 1\}$ as the sign function. By $\sigma \sim+$ (resp. $\sigma \sim-$), we say that the sign function $\sigma$ is equivalent to  all positive sign function (resp. all negative sign function). In a signed graph, the sign of a cycle  is defined to be the product of the signs of its edges. A signed cycle is said to be positive (resp. negative) if its sign is positive (resp. negative). A signed graph is said to be balanced if none of its cycles is negative, otherwise unbalanced.\\
 \indent In a graph $G$ with vertex set $V_G=\left\{v_{1}, \ldots, v_{n}\right\}$, if two distinct vertices $v_i$ and $ v_j$ are adjacent, we write $v_i \sim v_j$, otherwise, $v_i \nsim v_j$. The adjacency matrix $A_S$ of a signed graph $S$  is a  $(0,1,-1)$ square symmetric matrix of order $n$ in which $a_{i j}=\sigma(v_iv_j)$ if and only if $v_{i} \sim v_{j}$, where $\sigma(v_iv_j)$ is the sign of an edge $v_iv_j$.\\
   \indent Many well-known graph ideas are directly applicable to the domain of signed graphs. For example, a signed graph  is regular if the underlying graph is regular. Similarly, the degree of a vertex $v$ in $S$ is its degree in $G$. However, certain concepts are only applicable to signed graphs. The positive degree $d_{S}^{+}(v)$  of a vertex $v$ in $S$ is the number of positive edges incident on $v$. Similarly, we define the negative degree $d_{S}^{-}(v)$. In $S$, the net-degree of $v$  is defined as $d_{S}^{\pm}(v)=d_{S}^{+}(v)-d_{S}^{-}(v)$. A signed graph is said to be net-regular if the net-degree, considered as a function on the vertex set, is constant.  The Laplacian matrix of a signed graph $S$ is defined to be $ L_S  = D_S - A_S$, where $A_S$ and $D_S$ are, respectively, the adjacency matrix and the diagonal matrix of vertex degrees of $S$.   The net Laplacian matrix of a signed graph $S$ is defined as $N_{S}=D_{{S}}^{\pm}-A_{{S}}$, where $D_{{S}}^{\pm}$ is the diagonal matrix of net-degrees of $S$. The term ``net Laplacian'' appeared for the first time in \cite{Netz1}. Clearly, the net Laplacian coincides with the Laplacian in the case of unsigned graphs. \\
\indent For a square matrix $C $ having order $n$, denote by
$
\psi_C ( t)=\operatorname{det}\left(t I_{n}-C\right),
$
 the characteristic polynomial of $C$, where $I_{n}$ is the identity matrix of size $n$. For a signed graph $S$, we call $\psi_{A_S}(t)$  (resp., $\psi_{L_S}(t)$, $\psi_{N_S}(t)$) the adjacency (resp., Laplacian, net Laplacian) characteristic polynomial of $S$, and its roots the adjacency (resp., Laplacian, net Laplacian) eigenvalues of $S$. The collection of eigenvalues of $A_S$ together with their multiplicities is called the adjacency spectrum of $S$.  Similar terminology will be used for $L_S$ and $N_S$. If $X$ is a subset of the vertex set of $S$, then we denote by $S^{X}$,  the signed graph obtained by reversing the sign of every edge with exactly one end in $X$. We say that $S^{X}$ is switching equivalent to $S$.  Switching equivalent signed graphs have the same spectrum as the adjacency matrix and the Laplacian matrix. Interestingly, this does not hold true for the spectrum of the net Laplacian matrix. If the eigenvalues of  $C_S$ are written as $c_{1}(S) \leqslant c_{2}(S) \leqslant \cdots \leqslant c_{n}(S)$, then its spectrum will be denoted by $\sigma_C(S)= \{c_{1}(S), c_{2}(S),\ldots, c_{n}(S)\}.$ \\
 \indent It is well established that if a connected unsigned graph has only two distinct adjacency eigenvalues, then it must be a complete graph. Till now,  to characterize all unsigned graphs with three distinct adjacency eigenvalues is still an unsolved problem. Unlike unsigned graphs, if a signed graph  has only
two distinct adjacency eigenvalues, then it need not be a complete graph. Ramezami \cite{10} proved
that if a signed graph has two distinct adjacency eigenvalues, then it must be  regular. Signed graphs which are $2,3,4,5-$regular and have only
two distinct adjacency eigenvalues are characterized in \cite{e5,y7,z10,z11}. Signed graphs with 3 distinct adjacency eigenvalues may or may not be regular. Some results concerning signed graphs with 3 distinct adjacency eigenvalues can be found in \cite{a1,F1,F2}. Thus, it will be interesting to  find the signed graphs with $4$ and $5$ distinct adjacency eigenvalues  from the known signed graphs with $2$ and $3$ distinct adjacency eigenvalues. \\
\indent If two signed graphs have the same adjacency (resp., Laplacian, net Laplacian) spectrum, then they are said to be adjacency (resp., Laplacian, net Laplacian) cospectral. Any two switching isomorphic signed graphs are adjacency (resp., Laplacian) cospectral. A signed graph is said to be determined by its  adjacency spectrum (resp., Laplacian spectrum) if adjacency cospectral (resp., Laplacian cospectral) signed graphs are switching isomorphic signed graphs. It is said to be determined by its net Laplacian spectrum  if net Laplacian cospectral  signed graphs are isomorphic signed graphs. In general, the spectrum does not determine the signed graph and this
problem has pushed a lot of research. Thus, it will be interesting to identify
adjacency (resp., Laplacian, net Laplacian) cospectral  non-isomorphic signed graphs for a given class of signed graphs. \\
\indent The rest of the paper is organised as follows. Section $2$  deals with the definitions and adjacency spectra of s-neighbourhood corona of  two signed graphs. In Section $3$, we obtain the Laplacian spectra of s-neighbourhood corona of  two signed graphs. In Section $4$, we obtain the net Laplacian spectra of s-neighbourhood corona of  two signed graphs. In Section $5$, we find the signed graphs with $4$ and $5$ distinct adjacency (resp. Laplacian, net Laplacian) eigenvalues  from the known signed graphs with $2$ and $3$ distinct adjacency (resp. Laplacian, net Laplacian) eigenvalues.
 \section{Adjacency spectra  of s-neighbourhood corona of two  signed graphs}
 In this section, we determine the adjacency spectra  of s-neighbourhood corona of two  signed graphs.  For that we need the following motivation.\\
 \indent Let $A$ be a  square matrix of order $n$. The $A$-coronal, denoted by $\varkappa_{A}(t)$, of  a matrix $A$ is defined  to be the sum of the entries of the matrix $\left(t I_{n}-A\right)^{-1}$, that is,
$$
\varkappa_{A}(t)={j}_{n}^{T}\left(t I_{n}-A\right)^{-1} {j}_{n},
$$
where $j_{n}$ denotes the column vector of size $n$ with all the entries equal to one. For the detailed information about the coronal of a matrix, we refer the reader to \cite{cu, mc}.\\
\indent If $A$ is an $n \times n$ matrix with each row sum equal to a fixed real number $k$, then it is well-known \cite{mc} that
\begin{equation}\label{eq2.1}
\varkappa_{A}(t)=\frac{n}{t-k}.
\end{equation}
In particular, since for any signed graph $S$ with $n$ vertices, each row sum of net Laplacian matrix $N_S$ is equal to 0, therefore, we obtain
\begin{equation}\label{2.2}
\varkappa_{N_S}(t)=\frac{n}{t}.
\end{equation}
 \indent For two matrices $C=\left(c_{i j}\right)_{m \times n}$ and $D=\left(d_{i j}\right)_{p \times q}$, the Kronecker product $C \otimes D$ is the $m p \times n q$ matrix obtained from $C$ by replacing each element $c_{i j}$ by $c_{i j} D$. This is an associative operation with the property that $(C \otimes D)^{T}=C^{T} \otimes D^{T}$ and $(C \otimes D)(E \otimes F)=C E \otimes D F$ whenever the products $C E$ and $D F$ exist. The latter means that $(C \otimes D)^{-1}=C^{-1} \otimes D^{-1}$ for nonsingular matrices $C$ and $D$. Also, if $C$ and $D$ are $n \times n$ and $p \times p$ matrices, then $\operatorname{det}(C \otimes D)=(\operatorname{det} C)^{p}(\operatorname{det} D)^{n}$ and the Kronecker sum is $D \oplus C = C \otimes I_p + I_n \otimes D$, where $I_p$ and $I_n$ are identity matrices of order $p$ and $n$, respectively.\\
 \indent By Schur complement formula, the determinant  of a $2\times 2$ block matrix  is given  by
$$
\operatorname{det}\left(\begin{array}{ll}
X_1 & X_2 \\
X_3 & X_4
\end{array}\right)=\operatorname{det}(X_4)\operatorname{det}( X_1-X_2 X_4^{-1} X_3),
$$
where $X_1$ and $X_4$ are square blocks and $X_4$ is nonsingular. Moreover, $X_1-X_2 X_4^{-1} X_3$ is known as Schur complement of  $X_4$.\\
\indent Neighbourhood corona of two unsigned graphs  was introduced in \cite{Gop}. In the same paper, the adjacency spectra and Laplacian spectra of neighbourhood corona of any two unsigned graphs is expressed by the corresponding spectra of  two factor unsigned graphs. Here, we introduce the following definition for neighbourhood corona of two signed graphs.
\begin{definition}
 Let $S_{1}$ and $S_{2}$ be two signed graphs on $n_{1}$ and $n_{2}$ vertices, respectively.  The signed neighbourhood corona   of $S_{1}$ and $S_{2}$, denoted by $S_{1} \star_s S_{2}$, is the signed graph obtained by taking one copy of $S_{1}$ and $n_{1}$ copies of $S_{2}$  and joining every neighbour of the $i$th vertex of $S_{1}$  with the same sign as the sign of incident edge to every vertex in the $i$th copy of $S_{2}$. For illustration, see Figure $1$, where dotted lines are negative edges and bold lines are positive edges.
\end{definition}
It is worth noting that the sign of added edges in  $S_{1} \star_s S_{2}$ is defined by the sign function of the signed graph $ S_{1}=(G_1,\sigma)$. This inspires us to call it signed neighbourhood corona (in short s-neighbourhood corona). Clearly, the signed graph $S_{1} \star_s S_{2}$  has $n_1(n_2+1)$ vertices. \\
\indent The following result will be useful to determine the adjacency spectrum of the signed graph $S_{1} \star_s S_{2}$.
\begin{theorem}\label{t2.2}
Let $S_{1} \star_s S_{2}$ be s-neighbourhood corona of  signed graphs $S_{1}$ and $S_{2}$ with  $n_1$ and $n_2$ vertices, respectively. If $S=S_{1} \star_s S_{2}$, then
$$ \psi_{A_S}(t)=\left(\psi_{A_{S_{2}}}(t)\right)^{n_{1}} \psi_{(A_{S_{1}}+\varkappa_{A_{S_{2}}}(t)A_{S_{1}}^2)}(t), $$
where $\varkappa_{A_{S_{2}}}(t)$ is the coronal of  adjacency matrix of a signed graph $S_2$.
\end{theorem}
\noindent {\bf Proof.} Let $S_{1}$ and $S_{2}$ be  signed graphs on $n_{1}$ and $n_{2}$ vertices, respectively. We first label the vertices of $S=S_{1} \star_s S_{2}$ as follows. Let $V_{S_{1}}=\left\{v_{1}, v_{2}, \ldots, v_{n_{1}}\right\}$ and $V_{S_{2}}=\left\{u_{1}, u_{2}, \ldots, u_{n_{2}}\right\}$. Let $\{u_{1}^{j}, u_{2}^{j}, \ldots, u_{n_{2}}^{j}\}$, $j=1,2, \ldots, n_{1}$, denote the vertex set of the $j$th copy of $S_{2}$, with the understanding that $u_{i}^{j}$ is the copy of $u_{i}$ for each $i$. Denote by
$$
V_{i}=\left\{u_{i}^{1}, u_{i}^{2}, \ldots, u_{i}^{n_{1}}\right\}, \quad i=1,2, \ldots, n_{2}.
$$
Then $V_{S_{1}} \cup V_{1} \cup V_{2} \cup \cdots \cup V_{n_{2}}$ is a partition of $V_S$. With this partition, the adjacency matrix of the signed graph $S_{1} \star_s S_{2}$ is given by
$$
A_S =\left(\begin{array}{cc}
A_{S_1} & j_{n_{2}}^T \otimes A_{S_1} \\
\left(j_{n_{2}}^T \otimes A_{S_1}\right)^{T} & A_{S_2} \otimes I_{n_{1}}
\end{array}\right).
$$
\begin{figure}
\centering
	\includegraphics[scale=1.2]{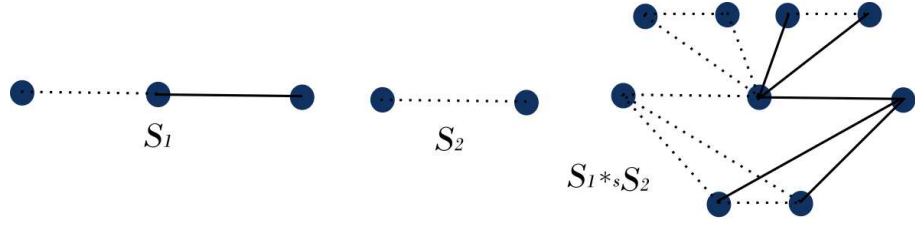}
	\caption{Signed nieghbourhood corona of two signed graphs}
	\label{biFigure6}
\end{figure}
By Schur complement formula and  viewed as a matrix over the field of rational functions $C(t)$, we obtain
\begin{equation*}\label{2.3}
\begin{aligned}
\psi_{A_S}(t)&= \operatorname{det}\left(t I_{n_{1}\left(n_{2}+1\right)}-A_S\right) \\
&= \operatorname{det}\left(\begin{array}{cc}
t I_{n_{1}}-A_{S_1} &  -(j_{n_{2}}^{T} \otimes A_{S_1}) \\
-\left(j_{n_{2}}^{T} \otimes A_{S_1}\right)^{T} & t I_{n_{1} n_{2}}-A_{S_{2}} \otimes I_{n_{1}}
\end{array}\right) \\
&= \operatorname{det}\left(\begin{array}{cc}
t I_{n_{1}}-A_{S_1} & - (j_{n_{2}}^{T} \otimes A_{S_1}) \\
-\left(j_{n_{2}}^{T} \otimes A_{S_1}\right)^{T} & (t I_{ n_{2}}-A_{S_{2}}) \otimes I_{n_{1}}
\end{array}\right)\\
&=\operatorname{det}\left(\left(t I_{n_{2}}-A_{S_{2}}\right) \otimes I_{n_{1}}\right) \times \operatorname{det}\{\left.t I_{n_{1}}-A_{S_1}\right.
-\left(j_{n_{2}}^{T} \otimes A_{S_1}\right)\left((t I_{ n_{2}}-A_{S_{2}}) \otimes I_{n_{1}}\right)^{-1}\left(j_{n_{2}} \otimes A_{S_1}\right)\}\\
&=\left(\operatorname{det}\left(t I_{n_{2}}-A_{S_{2}}\right)\right)^{n_{1}} \times \operatorname{det}\left \{\left.t I_{n_{1}}-A_{S_1}\right.\left.\left.-\left(j_{n_{2}}^{T}\left(t  I_{n_{2}}-A_{S_{2}}\right)^{-1}\right. j_{n_{2}}\right. \otimes A_{S_1}^2\right)\right \} \\
&= \left(\operatorname{det}\left(t I_{n_{2}}-A_{S_{2}}\right)\right)^{n_{1}} \operatorname{det}\left(t I_{n_{1}}-A_{S_{1}}-\varkappa_{A_{S_{2}}}(t)A_{S_{1}}^2\right)\\
&=\left(\psi_{A_{S_{2}}}(t)\right)^{n_{1}} \psi_{(A_{S_{1}}+\varkappa_{A_{S_{2}}}(t)A_{S_{1}}^2)}(t).
\end{aligned}
\end{equation*}
Hence the result follows. \qed
 \begin{theorem}\label{t2.3}
 Let $S_{1}$ be a signed graph on $n_{1}$ vertices and $S_{2}$ be an $r_{2}-$net regular signed graph on $n_{2}$ vertices. Assume that $\sigma_A\left(S_{1}\right)=\left\{\theta_{1}(S_1), \theta_{2}(S_1), \ldots, \theta_{n_{1}}(S_1)\right\}$ and $\sigma_A\left(S_{2}\right)=$ \\ $\left\{\theta_{1}(S_2), \theta_{2}(S_2), \ldots, \theta_{n_{2}}(S_2)\right\}$, where $\theta_{k}(S_2)=r_2$ for fixed integer $k$, $1\leq k \leq n_2$.  Then, the adjacency spectrum of $S_{1} \star_s S_{2}$ is given by\\
$(i)$ the adjacency eigenvalue $\theta_{j}(S_2)$ with multiplicity $n_{1}$ for every adjacency eigenvalue $\theta_{j}(S_2)(j=1,2, \ldots, n_{2}$ and $j\neq k)$ of $\sigma_A\left(S_{2}\right)$,\\
$(ii)$ two adjacency eigenvalues
$$
\frac{\theta_{i}\left(S_{1}\right)+r_{2} \pm \sqrt{(1+4 n_{2}) \theta_{i}\left(S_{1}\right)^{2}+\left(r_2-2\theta_{i}\left(S_{1}\right)\right)r_{2}}}{2}
$$
 for each adjacency eigenvalue $\theta_{i}(S_1) \left(i=1,2, \ldots, n_{1}\right)$ of $\sigma_A\left(S_{1}\right)$.
 \end{theorem}

\noindent {\bf Proof.} By Theorem \ref{t2.2}, we have
\begin{equation}\label{eqt2.3}
\begin{aligned} \psi_{A_S}(t)&=\left(\psi_{A_{S_{2}}}(t)\right)^{n_{1}} \psi_{(A_{S_{1}}+\varkappa_{A_{S_{2}}}(t)A_{S_{1}}^2)}(t)\\
&= \left(\operatorname{det}\left(t I_{n_{2}}-A_{S_{2}}\right)\right)^{n_{1}} \operatorname{det}\left(t I_{n_{1}}-A_{S_{1}}-\varkappa_{A_{S_{2}}}(t)A_{S_{1}}^2\right).
\end{aligned}
 \end{equation}
As $S_{2}$ is an $r_{2}-$net regular signed graph with $n_{2}$ vertices, therefore by Eq. (\ref{eq2.1}), we get
$$
\varkappa_{A_{S_{2}}}(t)=\frac{n_{2}}{t-r_{2}}.
$$
Clearly, the only pole of $\varkappa_{A_{S_{2}}}(t)$ is $t=r_{2}$. Thus, by Eq. (\ref{eqt2.3}), $\theta_{j}(S_2)$ is the adjacency eigenvalue with multiplicity $n_{1}$ for every adjacency eigenvalue $\theta_{j}(S_2)(j=1,2, \ldots, n_{2}$ and $j\neq k)$ of $\sigma_A\left(S_{2}\right)$. Now, the remaining $2 n_{1}$ adjacency eigenvalues of $S_{1} \star_s S_{2}$ are obtained by solving
$$
t-\theta_{i}\left(S_{1}\right)-\frac{n_{2}}{t-r_{2}} \theta_{i}\left(S_{1}\right)^{2}=0
$$
for each adjacency eigenvalue $\theta_{i}(S_1) \left(i=1,2, \ldots, n_{1}\right)$ of $\sigma_A\left(S_{1}\right)$. Hence  the result is proved. \qed
\begin{theorem}\label{c2.4} Let $S$ be a signed graph on $n$ vertices. For positive integers $p$ and $q$, let $(K_{p, q},-)$ denotes the signed graph with all negative function, where $K_{p, q}$ is a complete bipartite graph on $p+q$ vertices. Assume that $\sigma_A\left(S\right)=\left\{\theta_{1}(S), \theta_{2}(S), \ldots, \theta_{n}(S)\right\}$. Then the adjacency spectrum of $S \star_s (K_{p, q},-)$ is given by\\
$(i)$  $0$ with multiplicity $n(p+q-2)$,\\
$(ii)$ the three roots of the equation
$
t^{3}-\theta_{i}(S) t^{2}-\left(p q+(p+q) \theta_{i}(S)^{2}\right) t+p q \theta_{i}(S)\left(2 \theta_{i}(S)-1\right)=0
$
for each adjacency eigenvalue $\theta_{i}(S) \left(i=1,2, \ldots, n\right)$ of $\sigma_A\left(S\right)$.
\end{theorem}
\noindent {\bf Proof.} By Theorem \ref{t2.2}, we have
\begin{equation}\label{eqt2.4}
\begin{aligned} \psi_{A_{S \star_s (K_{p, q},-)} }(t)= \left(\operatorname{det}\left(t I_{p+q}-A_{(K_{p, q},-)}\right)\right)^{n} \operatorname{det}\left(t I_{n}-A_{S}-\varkappa_{A_{(K_{p, q},-)}}(t)A_{S}^2\right).
\end{aligned}
 \end{equation}
 With a suitable labelling of vertices, the adjacency matrix of a signed graph $(K_{p, q},-)$ is given by $$A_{(K_{p, q},-)}=\left(\begin{array}{ll}{O}_{ p} & -{j}_{p\times q} \\ -{j}_{q\times p} & {O}_{ q}\end{array}\right),$$
 where $O_p$ is a zero matrix of order $p$ and $j_{p\times q}$ is a matrix of order $p\times q$ with each entry equal to one.\\  Now, let $Z=\operatorname{diag}\left((t-q) I_p,(t-p) I_q\right)$ be the diagonal matrix with  first $p$ diagonal entries being $(t-q)$ and the last $q$ entries being $(t-p)$. Then $(t I_{p+q}-A_{(K_{p, q},-)}) Z j_{p+q}=\left(t^2-p q\right) j_{p+q}$ and so
$$
\varkappa_{A_{(K_{p, q},-)}}(t)=j_{p+q}^T(t I_{p+q}-A_{(K_{p, q},-)})^{-1} j_{p+q}=\frac{j_{p+q}^T Z j_{p+q}}{t^2-p q}=\frac{(p+q)t-2 p q}{t^2-p q}.
$$
The adjacecny spectrum \cite{ts} of $(K_{p, q},-)$ consists of eigenvalues $\pm \sqrt{p q}$ with multiplicity $1$, and $0$ with multiplicity $p+q-2$. Since, the poles of $\varkappa_{A_{\left(K_{p, q},-\right)}}(t)$ are $\pm \sqrt{p q}$, hence  the result follows by Eq. (\ref{eqt2.4}). \qed
\indent Proceeding similarly as in Theorem \ref{c2.4}, we have the following observation.
\begin{theorem} Let $S$ be a signed graph on $n$ vertices. For positive integers $p$ and $q$, let $(K_{p, q},+)$ denotes the signed graph with all positive function, where $K_{p, q}$ is a complete bipartite graph on $p+q$ vertices. Assume that $\sigma_A\left(S\right)=\left\{\theta_{1}(S), \theta_{2}(S), \ldots, \theta_{n}(S)\right\}$.  Then the adjacency spectrum of $S \star_s (K_{p, q},+)$ is given by\\
$(i)$  $0$ with multiplicity $n(p+q-2)$,\\
$(ii)$ the three roots of the equation
$
t^{3}-\theta_{i}(S) t^{2}-\left(p q+(p+q) \theta_{i}(S)^{2}\right) t-p q \theta_{i}(S)\left(2 \theta_{i}(S)-1\right)=0
$
for each adjacency eigenvalue $\theta_{i}(S) \left(i=1,2, \ldots, n\right)$ of $\sigma_A\left(S\right)$.
\end{theorem}
\section{Laplacian spectra  of s-neighbourhood corona of two  signed graphs}
This section begins with the following main result.
 \begin{theorem}\label{t3.1}
Let $S_{1}$ and $S_{2}$ be two signed graphs on $n_{1}$ and $n_{2}$ vertices, respectively. Also, let $P$ be a nonsingular matrix such that $
t I_{ n_{2}}-L_{S_{2}}=P \Omega P^{-1},$
where $\Omega$ is a diagonal matrix whose diagonal entries are the Laplacian eigenvalues of $S_2$. If $S=S_{1} \star_s S_{2}$, then
$$
\begin{array}{rlr}
\psi_{L_S}(t)  =\operatorname{det}\left(- D_{{S_1}} \oplus (tI_{n_2}-L_{S_{2}}) \right) \times \operatorname{det}(M),
\end{array}
$$
where $\operatorname{det} (M)
=$
\begin{equation*}
\begin{aligned}
 \operatorname{det}\left\{tI_{n_1}-L_{S_{1}}- n_{2} D_{{S_1}}-\left(j_{n_{2}}^{T}P\otimes -A_{S_{1}}D_{{S_1}}^{-1}\right)
 {\left(\Omega \otimes - D_{{S_1}}^{-1}+I_{n_2} \otimes I_{n_1}\right)^{-1}\left(P^{-1}j_{n_{2}}\otimes A_{S_{1}}\right)}\right \}.
\end{aligned}
\end{equation*}

\end{theorem}
\noindent {\bf Proof.} Using the  labelling as in Theorem \ref{t2.2}, the Laplacian matrix of signed graph $S=S_{1} \star_s S_{2}$ is given by
$$
L_{S}=\left(\begin{array}{ll}
L_{S_{1}}+ n_{2} D_{{S_1}} & -(j_{n_{2}}^{T} \otimes A_{S_{1}}) \\
-\left(j_{n_{2}}^{T} \otimes A_{S_{1}}\right)^{T} &   D_{{S_1}}\oplus L_{S_{2}}
\end{array}\right).
$$
 Therefore
$$
\begin{aligned}
\psi_{L_{S}}(t) &=\operatorname{det}\left(\begin{array}{ll}
tI_{n_1}-L_{S_{1}}- n_{2} D_{{S_1}} & (j_{n_{2}}^{T} \otimes A_{S_{1}}) \\
\left(j_{n_{2}}^{T} \otimes A_{S_{1}}\right)^{T} &   - D_{{S_1}} \oplus (tI_{n_2}-L_{S_{2}})
\end{array}\right) \\
&=\operatorname{det}\left(- D_{{S_1}} \oplus (tI_{n_2}-L_{S_{2}})\right) \times  \operatorname{det}(M),
\end{aligned}
$$
where  $M=\left(t I_{n_{1}+n_{1} n_{2}}-L_{S_{1} \star_s S_{2}}\right) /\left(-D_{S_1}\oplus(t I_{ n_{2}}-L_{S_{2}})  \right)$ is the Schur complement with respect to $\left( -D_{S_1}\oplus(t I_{ n_{2}}-L_{S_{2}}) \right)$.
Now,
\begin{equation}\label{ee2.5}
\begin{aligned}
\operatorname{det} (M)=& \operatorname{det}\left(\left(t I_{n_{1}+n_{1} n_{2}}-L_{S_{1} \star_s S_{2}}\right) /(-D_{S_1}\oplus(t I_{ n_{2}}-L_{S_{2}}))\right)\\
=& \operatorname{det}\left\{tI_{n_1}-L_{S_{1}}- n_{2} D_{{S_1}}-\left(j_{n_{2}}^{T}\otimes A_{S_1}\right)\right.
 {( -D_{S_1}\oplus(t I_{ n_{2}}-L_{S_{2}})  ) ^{-1}\left(j_{n_{2}}^{T}\otimes A_{S_1}\right)^{T}} \}.
\end{aligned}
\end{equation}
 Also,
\begin{equation}\label{inverse3.6}
\begin{aligned}
( -D_{S_1}\oplus(t I_{ n_{2}}-L_{S_{2}})  ) ^{-1}&=(t I_{ n_{2}}-L_{S_{2}} \otimes I_{n_1}+I_{n_2} \otimes -D_{S_1})^{-1}\\ &=\left(\left(t I_{ n_{2}}-L_{S_{2}}  \otimes -D_{S_1}^{-1}+I_{n_2} \otimes I_{n_1}\right)(I_{n_2} \otimes -D_{S_1})\right)^{-1} \\
&=\left(I_{n_2} \otimes -D_{S_1}^{-1}\right)\left(t I_{ n_{2}}-L_{S_{2}} \otimes -D_{S_1}^{-1}+I_{n_2} \otimes I_{n_1}\right)^{-1}.
\end{aligned}
\end{equation}
Eigen-decomposing $t I_{ n_{2}}-L_{S_{2}}$ through nonsingular matrix $P$. We have
\begin{equation}\label{eigen2.7}
t I_{ n_{2}}-L_{S_{2}}=P \Omega P^{-1},
\end{equation}
where $\Omega$ is the diagonal matrix whose diagonal entries are the Laplacian eigenvalues of $S_2$.
From Eq.s  (\ref{inverse3.6}) and (\ref{eigen2.7}), we have
\begin{equation}\label{cb2.8}
\begin{aligned}
(-D_{S_1}\oplus(t I_{ n_{2}}-L_{S_{2}})  )^{-1} &=\left(I_{n_2} \otimes -D_{S_1}^{-1}\right)\left(P \Omega P^{-1} \otimes  -D_{S_1}^{-1} +I_{n_2} \otimes I_{n_1}\right)^{-1} \\
&=\left(I_{n_2} \otimes -D_{S_1}^{-1}\right)\left(( P \otimes I_{n_1} ) ( \Omega  \otimes -D_{S_1} ^ { - 1 } + I_{n_2} \otimes I_{n_1} ) \left(P^{-1}\otimes I_{n_1}\right)\right)^{-1}\\
&=\left( P \otimes -D_{S_1}^{-1} \right)\left(\Omega \otimes -D_{S_1}^{-1}+I_{n_2} \otimes I_{n_1}\right)^{-1}\left(P^{-1} \otimes I_{n_1}\right).
\end{aligned}
\end{equation}
Substituting Eq. (\ref{cb2.8}) in Eq. (\ref{ee2.5}), we obtain  $\operatorname{det} (M)
=$
\begin{equation*}
\begin{aligned}
 \operatorname{det}\left\{tI_{n_1}-L_{S_{1}}- n_{2} D_{{S_1}}-\left(j_{n_{2}}^{T}P\otimes -A_{S_{1}}D_{{S_1}}^{-1}\right)
 {\left(\Omega \otimes - D_{{S_1}}^{-1}+I_{n_2} \otimes I_{n_1}\right)^{-1}\left(P^{-1}j_{n_{2}}\otimes A_{S_{1}}\right)}\right \}.
\end{aligned}
\end{equation*}
Hence the result follows. \qed
 \indent The following result gives the eigenvalues of Kronecker sum in terms of the eigenvalues of its factor matrices.
\begin{lemma}\label{l2.8}
\emph {\cite{Kron}} Let $A$ and $B$ be two square matrices of order $n$ and $m$, respectively. Assume that $x_i$, $i=1,2,\ldots,n$ and $y_j$, $j=1,2,\ldots,m$ are  eigenvalues of the matrices $A$ and $B$, respectively. Then the eigenvalues of the
Kronecker sum  $B \oplus A$ are  $x_i+y_j$, $i=1,2,\ldots, n$ and $ j=1, 2, \ldots, m$.

\end{lemma}
Next result completely determines   the Laplacian spectrum of $S_{1} \star_s S_{2}$  for regular $S_{1} $ and regular and net regular $ S_{2}$.
\begin{theorem}\label{t3.3}
 Let $S_{1}$ be an $r_1-$regular signed graph on $n_{1}$ vertices and $S_{2}$ be an $r_{2}-$ regular and $r_{3}-$net regular signed graph on $n_{2}$ vertices. Assume that $\sigma_L\left(S_{1}\right)=\left\{\lambda_{1}(S_1), \lambda_{2}(S_1), \ldots, \lambda_{n_{1}}(S_1)\right\}$ and $\sigma_L\left(S_{2}\right)=$  $\left\{\lambda_{1}(S_2), \lambda_{2}(S_2), \ldots, \lambda_{n_{2}}(S_2)\right\}$, where $\lambda_{k}(S_2)=r_2-r_3$ for fixed integer $k$, $1\leq k \leq n_2$.  Then, the Laplacian spectrum of $S_{1} \star_s S_{2}$ is given by\\
$(i)$ the Laplacian eigenvalue $\lambda_{j}(S_2)+r_1$ with multiplicity $n_{1}$ for every Laplacian eigenvalue $\lambda_{j}(S_2)(j=1,2, \ldots, n_{2}$ and $j\neq k)$ of $\sigma_L\left(S_{2}\right)$,\\
$(ii)$ two Laplacian eigenvalues
$$
\frac{\beta \pm \sqrt{\beta^{2}+4 \left( n_{2} (\lambda_{i}\left(S_{1}\right)-r_1)^2-(\lambda_{i}(S_{1})+r_1n_2)(r_1+r_2-r_3)\right)}}{2},
$$
 where $\beta = r_1+r_2-r_3 +\lambda_{i}(S_{1})+r_1n_2$,       for each Laplacian eigenvalue $\lambda_{i}(S_1) \left(i=1,2, \ldots, n_{1}\right)$ of $\sigma_L\left(S_{1}\right)$.
 \end{theorem}
\noindent {\bf Proof.} By Theorem \ref{t3.1}, we have
$$
\begin{array}{rlr}
\psi_{L_S}(t)  =\operatorname{det}\left(- D_{{S_1}} \oplus (tI_{n_2}-L_{S_{2}}) \right) \times \operatorname{det}(M),
\end{array}
$$
where $\operatorname{det} (M)
=$
\begin{equation}\label{e3.3}
\begin{aligned}
 \operatorname{det}\left\{tI_{n_1}-L_{S_{1}}- n_{2} D_{{S_1}}-\left(j_{n_{2}}^{T}P\otimes -A_{S_{1}}D_{{S_1}}^{-1}\right)\right.
 {\left(\Omega \otimes - D_{{S_1}}^{-1}+I_{n_2} \otimes I_{n_1}\right)^{-1}\left(P^{-1}j_{n_{2}}\otimes A_{S_{1}}\right)} \}.    
\end{aligned}
\end{equation}
 If $S_{1}$ is an $r_{1}-$regular signed graph, then we can  easily obtained  from  Eq. (\ref{e3.3}) and Lemma \ref{l2.8} that

$$
\begin{array}{rlr}
\psi_{L_{S_{1} \star_s S_{2}}}(t)  =(\operatorname{det}\left((t-r_1) I_{ n_{2}}-L_{S_{2}} \right))^{n_1} \times \operatorname{det}(M),
\end{array}
$$
where $\operatorname{det} (M)
=$
\begin{equation*}
\begin{aligned}
 &\operatorname{det}\left\{\left(t-r_{1} n_{2}\right) I_{n_{1}}-L_{S_{1}}-\left(j_{n_{2}}^{T}P\otimes -A_{S_1}(r_1I_{n_1})^{-1}\right)
 {\left((\frac{1}{r_1}\Omega-I_{n_2})   \otimes I_{m_1}\right)^{-1}\left(P^{-1}j_{n_{2}}\otimes A_{S_1}\right)} \right\}\\
 &= \operatorname{det}\{(t-r_{1} n_{2}) I_{n_{1}}-L_{S_{1}}-\varkappa_{L_{S_{2}}}(t-r_1)A_{S_1}^2
  \}\\
  &= \operatorname{det}\{(t-r_{1} n_{2}) I_{n_{1}}-L_{S_{1}}-\varkappa_{L_{S_{2}}}(t-r_1)(r_1I_{n_1}-L_{S_1})^2
  \}.
\end{aligned}
\end{equation*}
Thus,
$$
\begin{array}{rlr}
\psi_{L_{S_{1} \star_s S_{2}}}(t)  =(\operatorname{det}\left((t-r_1) I_{ n_{2}}-L_{S_{2}} \right))^{n_1} \times\operatorname{det}\{(t-r_{1} n_{2}) I_{n_{1}}-L_{S_{1}}-\varkappa_{L_{S_{2}}}(t-r_1)(r_1I_{n_1}-L_{S_1})^2
  \}.
\end{array}
$$
Since $S_{2}$ is an $r_{2}-$regular and $r_{3}-$net regular signed graph on $n_{2}$ vertices, therefore, the row sum of $L_{S_{2}}$ is  $r_{2}-r_3$. Also, by Eq. (\ref{eq2.1}), we get
$$ \varkappa_{L_{S_2}}(t)=\frac{n_2}{t-(r_{2}-r_3)}.$$
Clearly, the only pole of $\varkappa_{L_{S_{2}}}(t)$ is $t=r_{2}-r_3$. Thus $\lambda_{j}(S_2)+r_1$ is the Laplacian eigenvalue with multiplicity $n_{1}$ for every Laplacian eigenvalue $\lambda_{j}(S_2)(j=1,2, \ldots, n_{2}$ and $j\neq k)$ of $\sigma_L\left(S_{2}\right)$. Now, the remaining $2 n_{1}$ Laplacian eigenvalues of $S_{1} \star_s S_{2}$ are obtained by solving
$$
t-r_1n_2-\lambda_{i}\left(S_{1}\right)-\frac{n_{2}}{t-r_1-(r_{2}-r_3)} (r_1-\lambda_{i}\left(S_{1}\right))^{2}=0
$$
for each Laplacian eigenvalue $\lambda_{i}(S_1) \left(i=1,2, \ldots, n_{1}\right)$ of $\sigma_L\left(S_{1}\right)$. Hence  the result is proved. \qed
 If $S$ is  a connected signed graph with $n$ vertices, then it is well-known that $0$ is a Laplacian eigenvalue if and only if $S$ is balanced. Moreover, $0$ is a Laplacian eigenvalue corresponding to an eigenvector having each entry equal to one. Proceeding similarly as in Theorem \ref{t3.3}, we immediately obtain the following.
\begin{theorem}\label{t3.4}
 Let $S_{1}$ be an $r_1-$regular signed graph on $n_{1}$ vertices and $S_{2}$ be a connected balanced signed graph on $n_{2}$ vertices. Assume that $\sigma_L\left(S_{1}\right)=\left\{\lambda_{1}(S_1), \lambda_{2}(S_1), \ldots, \lambda_{n_{1}}(S_1)\right\}$ and $\sigma_L\left(S_{2}\right)=$  $\left\{\lambda_{1}(S_2), \lambda_{2}(S_2), \ldots, \lambda_{n_{2}}(S_2)\right\}$, where $\lambda_{k}(S_2)=0 $ for fixed integer $k$, $1\leq k \leq n_2$.  Then, the Laplacian spectrum of $S_{1} \star_s S_{2}$ is given by\\
$(i)$ the Laplacian eigenvalue $\lambda_{j}(S_2)+r_1$ with multiplicity $n_{1}$ for every Laplacian eigenvalue $\lambda_{j}(S_2)(j=1,2, \ldots, n_{2}$ and $j\neq k)$ of $\sigma_L\left(S_{2}\right)$,\\
$(ii)$ two Laplacian eigenvalues
$$
\frac{(n_2+1)r_1 +\lambda_{i}(S_{1}) \pm \sqrt{((n_2+1)r_1 +\lambda_{i}(S_{1}))^{2}-4\lambda_{i}(S_{1})\left((2n_2+1)r_1-n_2\lambda_{i}(S_{1}\right)}}{2},
$$
 for each Laplacian eigenvalue $\lambda_{i}(S_1) \left(i=1,2, \ldots, n_{1}\right)$ of $\sigma_L\left(S_{1}\right)$.
 \end{theorem}
Note that  Theorem \ref{t3.4} extends and generalizes the [Theorem $3.1$, \cite{Gop}] to signed graphs.
\section{Net Laplacian spectra  of s-neighbourhood corona of two  signed graphs}
This section begins with the  following result that will be useful to determine the net Laplacian spectrum of the signed graph $S_{1} \star_s S_{2}$.
\begin{theorem}\label{t4.1}
Let $S_{1}$ and $S_{2}$ be two signed graphs on $n_{1}$ and $n_{2}$ vertices, respectively. Also, let $P$ be a nonsingular matrix such that $
t I_{ n_{2}}-N_{S_{2}}=P \Omega P^{-1},$
where $\Omega$ is a diagonal matrix whose diagonal entries are the net Laplacian eigenvalues of $S_2$. If $S=S_{1} \star_s S_{2}$, then
$$
\begin{array}{rlr}
\psi_{N_S}(t)  =\operatorname{det}\left(- D^\pm_{{S_1}} \oplus (tI_{n_2}-N_{S_{2}}) \right) \times \operatorname{det}(M),
\end{array}
$$
where $\operatorname{det} (M)
=$
\begin{equation*}
\begin{aligned}
 \operatorname{det}\left\{tI_{n_1}-N_{S_{1}}- n_{2} D^\pm_{{S_1}}-\left(j_{n_{2}}^{T}P\otimes -A_{S_{1}}{D^\pm_{{S_1}}}^{-1}\right)
 {\left(\Omega \otimes - {D^\pm_{{S_1}}}^{-1}+I_{n_2} \otimes I_{n_1}\right)^{-1}\left(P^{-1}j_{n_{2}}\otimes A_{S_{1}}\right)} \right\}.
\end{aligned}
\end{equation*}

\end{theorem}
\noindent {\bf Proof.} Using the  labelling as in Theorem \ref{t2.2}, the net Laplacian matrix of the signed graph $S=S_{1} \star_s S_{2}$ is given by
$$
N_{S}=\left(\begin{array}{ll}
N_{S_{1}}+ n_{2} D^\pm_{{S_1}} & -(j_{n_{2}}^{T} \otimes A_{S_{1}}) \\
-\left(j_{n_{2}}^{T} \otimes A_{S_{1}}\right)^{T} &   D^\pm_{{S_1}}\oplus N_{S_{2}}
\end{array}\right).
$$
 Therefore
$$
\begin{aligned}
\psi_{N_{S}}(t) &=\operatorname{det}\left(\begin{array}{ll}
tI_{n_1}-N_{S_{1}}- n_{2} D^\pm_{{S_1}} & (j_{n_{2}}^{T} \otimes A_{S_{1}}) \\
\left(j_{n_{2}}^{T} \otimes A_{S_{1}}\right)^{T} &   - D^\pm_{{S_1}} \oplus (tI_{n_2}-N_{S_{2}})
\end{array}\right) \\
&=\operatorname{det}\left(- D^\pm_{{S_1}} \oplus (tI_{n_2}-N_{S_{2}})\right) \times  \operatorname{det}(M),
\end{aligned}
$$
where  $M=\left(t I_{n_{1}+n_{1} n_{2}}-N_{S_{1} \star_s S_{2}}\right) /\left(-D_{S_1}^\pm\oplus(t I_{ n_{2}}-N_{S_{2}})  \right)$ is the Schur complement with respect to $\left( -D_{S_1}^\pm\oplus(t I_{ n_{2}}-N_{S_{2}}) \right)$.
Now,
\begin{equation*}
\begin{aligned}
\operatorname{det} (M)=& \operatorname{det}\left(\left(t I_{n_{1}+n_{1} n_{2}}-N_{S_{1} \star_s S_{2}}\right) /(-D_{S_1}^\pm\oplus(t I_{ n_{2}}-N_{S_{2}}))\right)\\
=& \operatorname{det}\left\{tI_{n_1}-N_{S_{1}}- n_{2} D^\pm_{{S_1}}-\left(j_{n_{2}}^{T}\otimes A_{S_1}\right)
 {( -D_{S_1}^\pm\oplus(t I_{ n_{2}}-N_{S_{2}})  ) ^{-1}\left(j_{n_{2}}^{T}\otimes A_{S_1}\right)^{T}}\right \}.
\end{aligned}
\end{equation*}
 Also,
\begin{equation*}
\begin{aligned}
( -D_{S_1}^\pm\oplus(t I_{ n_{2}}-N_{S_{2}})  ) ^{-1}&=(t I_{ n_{2}}-N_{S_{2}} \otimes I_{n_1}+I_{n_2} \otimes -D_{S_1}^\pm)^{-1}\\ &=\left(\left(t I_{ n_{2}}-N_{S_{2}}  \otimes -{D^\pm_{S_1}}^{-1}+I_{n_2} \otimes I_{n_1}\right)(I_{n_2} \otimes -D_{S_1}^\pm)\right)^{-1} \\
&=\left(I_{n_2} \otimes -{D_{S_1}^\pm}^{-1}\right)\left(t I_{ n_{2}}-N_{S_{2}} \otimes -{D_{S_1}^\pm}^{-1}+I_{n_2} \otimes I_{n_1}\right)^{-1}.
\end{aligned}
\end{equation*}
Now, eigen-decomposing $t I_{ n_{2}}-N_{S_{2}}$ through nonsingular matrix $P$ such  that $
t I_{ n_{2}}-L_{S_{2}}=P \Omega P^{-1},$
where $\Omega$ is a diagonal matrix whose diagonal entries are the net Laplacian eigenvalues of $S_2$ and proceeding similarly as in Theorem \ref{t3.1}, we get
$\operatorname{det} (M)
=$
\begin{equation*}
\begin{aligned}
 \operatorname{det}\left\{tI_{n_1}-N_{S_{1}}- n_{2} D^\pm_{{S_1}}-\left(j_{n_{2}}^{T}P\otimes -A_{S_{1}}{D^\pm_{{S_1}}}^{-1}\right)
 {\left(\Omega \otimes - {D^\pm_{{S_1}}}^{-1}+I_{n_2} \otimes I_{n_1}\right)^{-1}\left(P^{-1}j_{n_{2}}\otimes A_{S_{1}}\right)} \right\}.
\end{aligned}
\end{equation*}
Hence the result follows. \qed

\begin{theorem}\label{t4.2}
Let $S_{1}$ be an  $r_2-$net regular signed graph with  $n_{1}$ vertices and $S_{2}$ be any signed graph with $n_{2}$ vertices.  Assume that $\sigma_N\left(S_{1}\right)=\left\{v_{1}(S_1), v_{2}(S_1), \ldots, v_{n_{1}}(S_1)\right\}$ and $\sigma_N\left(S_{2}\right)=$ $\left\{v_{1}(S_2), v_{2}(S_2), \ldots, v_{n_{2}}(S_2)\right\}$, where $v_{k}(S_2)=0$ for fixed integer $k$, $1\leq k \leq n_2$. Then the net Laplacian spectrum of $S_{1} \star_s S_{2}$ is given by\\
$(i)$ the net Laplacian eigenvalue $v_{j}(S_2)+ r_2$ with multiplicity $n_{1}$ for every net Laplacian eigenvalue $v_{j}(S_2)(j=1,2, \ldots, n_{2}$ and $j\neq k)$ of $\sigma_N\left(S_{2}\right)$,\\
$(ii)$ two multiplicity-one net Laplacian eigenvalues $$
\begin{aligned}
&\frac{v_{i}\left(S_{1}\right)+ \left(n_{2}+1\right) r_{2} \pm \sqrt{\left(v_{i}\left(S_{1}\right)+ \left(n_{2}+1\right) r_{2}\right)^{2}-4 v_{i}\left(S_{1}\right)\left(\left(2 n_{2}+1 \right) r_{2}-n_{2} v_{i}\left(S_{1}\right)\right)}}{2}
\end{aligned}
$$
 for each net Laplacian eigenvalue $v_{i}(S_1) \left(i=1,2, \ldots, n_{1}\right)$ of $\sigma_N\left(S_{1}\right)$.
\end{theorem}
\noindent {\bf Proof.} Let $S_{1}$ be an $r_{2}-$net regular signed graph. Then by Theorem \ref{t4.1} and Lemma \ref{l2.8}, we obtain

$$
\begin{array}{rlr}
\psi_{N_{S_{1} \star_s S_{2}}}(t)  =(\operatorname{det}\left((t-r_2) I_{ n_{2}}-N_{S_{2}} \right))^{n_1} \times \operatorname{det}(M),
\end{array}
$$
where $\operatorname{det} (M)
=$
\begin{equation*}
\begin{aligned}
 &\operatorname{det}\left\{\left(t-r_{2} n_{2}\right) I_{n_{1}}-N_{S_{1}}-\left(j_{n_{2}}^{T}P\otimes -A_{S_1}(r_2I_{n_1})^{-1}\right)
 {\left((\frac{1}{r_2}\Omega-I_{n_2})   \otimes I_{m_1}\right)^{-1}\left(P^{-1}j_{n_{2}}\otimes A_{S_1}\right)} \right\}\\
 &= \operatorname{det}\{(t-r_{2} n_{2}) I_{n_{1}}-N_{S_{1}}-\varkappa_{N_{S_{2}}}(t-r_2)A_{S_1}^2
  \}\\
  &= \operatorname{det}\{(t-r_{2} n_{2}) I_{n_{1}}-N_{S_{1}}-\varkappa_{N_{S_{2}}}(t-r_2)(r_2I_{n_1}-N_{S_1})^2
  \}.
\end{aligned}
\end{equation*}
Thus, we have
$$
\begin{array}{rlr}
\psi_{N_{S_{1} \star_s S_{2}}}(t)  =(\operatorname{det}\left((t-r_2) I_{ n_{2}}-N_{S_{2}} \right))^{n_1} \times\operatorname{det}\{(t-r_{2} n_{2}) I_{n_{1}}-N_{S_{1}}-\varkappa_{N_{S_{2}}}(t-r_2)(r_2I_{n_1}-N_{S_1})^2
  \}.
\end{array}
$$
Since  each row sum of the net Laplacian matrix $N_{S_2}$ is equal to 0, therefore, by Eq. (\ref{2.2}), we get
$$ \varkappa_{N_{S_2}}(t)=\frac{n_2}{t}.$$
Clearly, the only pole of $\varkappa_{N_{S_{2}}}(t)$ is $t=0$. Thus $v_{j}(S_2)+r_2$ is the net Naplacian eigenvalue with multiplicity $n_{1}$ for every net Naplacian eigenvalue $\lambda_{j}(S_2)(j=1,2, \ldots, n_{2}$ and $j\neq k)$ of $\sigma_N\left(S_{2}\right)$. Now, the remaining $2 n_{1}$ net Naplacian eigenvalues of $S_{1} \star_s S_{2}$ are obtained by solving
$$
t-r_2n_2-v_{i}\left(S_{1}\right)-\frac{n_{2}}{t-r_2} (r_2-\lambda_{i}\left(S_{1}\right))^{2}=0
$$
for each net Laplacian eigenvalue $\lambda_{i}(S_1) \left(i=1,2, \ldots, n_{1}\right)$ of $\sigma_N\left(S_{1}\right)$. Hence  the result is proved.\qed
\section{Construction of signed graphs with $4$ and $5$ distinct eigenvalues and spectral determination}
In unsigned graphs, there is a well-known relationship between the diameter and the number of distinct adjacency eigenvalues. The number of distinct adjacency eigenvalues cannot be less than the diameter plus 1. However, in signed graphs, in general,  this is not true. To see this, consider an unbalanced cycle on $4$ vertices. \\
\indent In the past $2-3$  years, there has been much investigation of signed graphs  with a few  distinct eigenvalues.  The s-neighbourhood corona construction can help us to utilize the known signed graphs  with $2$ and $3$ distinct eigenvalues to obtain new signed graphs  with $4$ and $5$  distinct eigenvalues. The following result is the direct consequence of Theorems \ref{t2.3}, \ref{t3.3} and \ref{t4.2}.
\begin{theorem}
$(i)$ Let $S_{1}$ be any signed graph on $n_1$ vertices with $t_1$ $(1\leq t_1\leq n_1)$ distinct adjacency eigenvalues and $S_{2}$ be a net regular signed graph on $n_2$ vertices  with $t_2$ $(1\leq t_2\leq n_2)$ distinct adjacency eigenvalues.  Then, the adjacency spectrum of $S_{1} \star_s S_{2}$ consists of atmost $2t_1+t_2$ distinct eigenvalues.\\
$(ii)$ Let $S_{1}$ be any regular signed graph on $n_1$ vertices with $t_1$ $(1\leq t_1\leq n_1)$ distinct Laplacian eigenvalues and $S_{2}$ be  regular and net regular signed graph on $n_2$ vertices  with $t_2$ $(1\leq t_2\leq n_2)$ distinct Laplacian eigenvalues.  Then, the Laplacian spectrum of $S_{1} \star_s S_{2}$ consists of atmost $2t_1+t_2$ distinct Laplacian eigenvalues.\\
$(iii)$ Let $S_{1}$ be any net regular signed graph on $n_1$ vertices with $t_1$ $(1\leq t_1\leq n_1)$ distinct net Laplacian eigenvalues and $S_{2}$ be any signed graph on $n_2$ vertices  with $t_2$ $(1\leq t_2\leq n_2)$ distinct net Laplacian eigenvalues.  Then, the  net Laplacian  spectrum of $S_{1} \star_s S_{2}$ consists of atmost $2t_1+t_2$ distinct  net Laplacian eigenvalues.
 \end{theorem}
 \begin{figure}
\centering
	\includegraphics[scale=1.2]{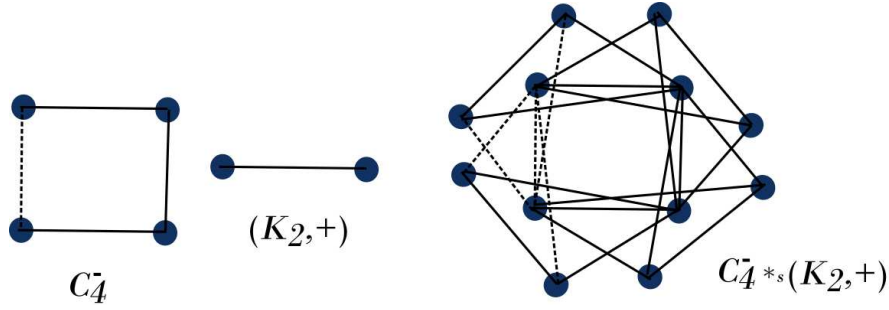}
	\caption{Signed graph $C_4^-\star_s (K_2,+)$ with a small number of adjacency eigenvalues. }
	\label{biFigure6}
\end{figure}
It is well established that if a signed graph on $n$ vertices has two distinct adjacency eigenvalues, then its spectrum is of the form $\{ \theta_1^a,\theta_2^b\}$, where $\theta_1\neq \theta_2$ and $a+b=n$. Thus, the following result directly follows from Theorem \ref{t2.3}.
\begin{theorem}\label{co}
$(i)$ Let $S$ be any signed graph on $n\geq 2$ vertices with $2$  distinct adjacency eigenvalues.   Then, the adjacency spectrum of $S \star_s (K_1,\sigma)$, where $K_1$ is a complete graph on $1$ vertex,  consists of $4$ distinct eigenvalues.\\
$(ii)$ Let $S$ be any signed graph on $n\geq 2$ vertices with $2$  distinct adjacency eigenvalues.   Then, the adjacency spectrum of $S \star_s (K_2,\sigma)$, where $K_2$ is a complete graph on $2$ vertices,  consists of $5$ distinct eigenvalues.
 \end{theorem}
\noindent {\bf Example.} Let $C_4^-$ be an unbalanced cycle on $4$ vertices. By Theorem \ref{co}, the signed graph $C_4^-\star_s (K_2,+)$ (shown in Figure $2$) on $12$ vertices has $5$ distint adjacency eigenvalues. Moreover, by Theorem \ref{t2.3}, the adjacency spectrum of $C_4^-\star_s (K_2,+)$ is given by $\left \{ -1^4,\frac{3\pm \sqrt{33}}{2}^2,\frac{-1\pm \sqrt{41}}{2}^2\right \}.$ \\
\indent Before to conclude this section, we will show that the signed graph $S_{1} \star_s S_{2}$ for arbitrary signed graphs $S_1$ and  $S_2$ is not determined by its adjacency, Laplacian and net Laplacian spectrum.  The following theorem is a direct consequence of Theorems \ref{t2.2}, \ref{t3.1} and \ref{t4.1}.
\begin{theorem}
$(i)$ Let $S_{1}$ and $S_{2}$  be two adjacency cospectral and  non-isomorphic signed graphs.   Then, the signed graphs $S_{1} \star_s S$ and $S_{2} \star_s S$ are non-isomorphic and adjacency cospectral for any arbitrary signed graph $S$.\\
$(ii)$ Let $S_{1}$ and $S_{2}$  be two Laplacian cospectral and non-isomorphic  signed graphs.   Then, the signed graphs $S_{1} \star_s S$ and $S_{2} \star_s S$ are non-isomorphic and Laplacian cospectral for any arbitrary signed graph $S$.\\
$(iii)$ Let $S_{1}$ and $S_{2}$  be two net Laplacian cospectral and non-isomorphic signed graphs.   Then, the signed graphs $S_{1} \star_s S$ and $S_{2} \star_s S$ are non-isomorphic and net Laplacian cospectral for any arbitrary signed graph $S$.\\
 \end{theorem}

\noindent{\bf Acknowledgements.}   This research is supported by SERB-DST research project number CRG/2020/000109. The research of Tahir Shamsher is supported by SRF financial assistance by Council of Scientific and Industrial Research (CSIR), New Delhi, India.

\noindent{\bf Conflict of interest.} The authors declare that they have no conflict of interest.\\

\noindent{\bf Data Availibility} Data sharing is not applicable to this article as no datasets were generated or analyzed
during the current study.\\

\end{document}